\documentclass[12pt]{article}
\usepackage{e-jc}
\usepackage{amsthm,
 }

\theoremstyle{plain}
\newtheorem{theorem}{Theorem}
\newtheorem{lemma}{Lemma}[section]
\newtheorem{corollary}[theorem]{Corollary}
\newtheorem{proposition}[theorem]{Proposition}

\theoremstyle{definition}

\theoremstyle{remark}

\def\c#1{\ensuremath{C_{#1}}}
\def\p#1{\ensuremath{P_{#1}}}
\def\3#1#2#3{_{#1=#2}^{#3}}   
\def\ca#1{\left|{#1}\right|}
\def\idn#1{\mathop\alpha(#1)}
\def\cqn#1{\mathop\omega(#1)}

\def\o#1#2{N({#1}|{#2})}
\def\mval#1 {\delta\mathopen({#1}\mathclose)}
\def\mxval#1 {\Delta\mathopen({#1}\mathclose)}
\def\dist{\mathop{\rm dist}\nolimits}
\def\B#1#2{\ensuremath{B({#2};{#1})}}
\def\girth#1{\mathop{\rm girth}({#1})}
\def\implies{\; \Longrightarrow \;}
\let\0=\emptyset
\def\N{{\bf N}}
\def\comp#1{\mathop{\rm comp}({#1})}
\def\lk{\mathop{{\rm lk}}\nolimits}
\let\mv=\delta
\def\crt{{\mathop{\rm cr}\nolimits}}
\def\e#1 {\ensuremath{\mathop{\rm
 e}\mathopen({#1}\mathclose)}}
\def\rr#1#2{\ensuremath{\mathop{\rm R}
 \mathopen({#1},{#2}\mathclose)}}
\let\ep=\varepsilon
\def\nv#1#2{\ensuremath{\##1(#2)}}

\def\DW#1{\ensuremath{BC_{#1}}}
\def\w#1{\ensuremath{Ch_{#1}}}
\def\5{\ensuremath{{\cal W}_{13;1,5}}}

\def\v#1{\ensuremath{d(#1)}}
\def\vind#1#2{\ensuremath{d_{#1}(#2)}}
\def\vv#1#2{\ensuremath{d^{\,#1}(#2)}}
\def\vvind#1#2#3{\ensuremath{d_{#1}^{#2}(#3)}}

\def\kdots{\ifmmode,\ldots,\else,~$\ldots\,$, \fi}
\def\row#1#2{\ifmmode #1_1\kdots #1_{#2}\else
 $#1_1$\kdots$#1_{#2}$\fi}

\def\noproof{\hfill$\square$}

\long\def\clause#1#2{\par\smallbreak\hangafter=1\hangindent
 20pt\noindent{\hbox to20pt{{#1\hfill}}{#2}}\hfill\par
 \ifdim \lastskip <\smallskipamount \removelastskip
 \penalty 15\smallskip \fi}

\def\loch{loop-chain}

\title{\bf Edge number critical triangle free
graphs with low independence numbers}

\author{J\"orgen Backelin}

\begin{document}
 \maketitle

 \begin{abstract}
 The structure of all triangle free graphs $G =
(V,E)$ with $\mathopen| E \mathclose| - 6 \mathopen|
V \mathclose| + 13\mathop\alpha (G) = 0$ is
determined, yielding an affirmative answer to a
question of Stanis\l aw Radziszowsky and Donald
Kreher.

\end{abstract}

 \section{Background.}
 \label{S:bg}
 A graph is {\sl edge number critical\/} (under
certain conditions), if it has the minimal possible
number of edges for all graphs fulfilling these
conditions. In this article, some edge number
critical graphs are determined among the triangle
free graphs (graphs without $K_3$ subgraphs) with a
prescribed number of vertices, upper bound for the
independence number, and sometimes also a prescribed
minimal valency.

 In my opinion, characterising such graphs has some
interest in itself. Moreover, often it is crucial
for determining better bounds for Ramsey numbers.

 In 1991, in \cite{RK}, Radziszowski and Kreher proved
that
 $$\ca E-6\ca V+13\idn G\ge0 \leqno(1)$$
 for any triangle free simple graph $G = (V,E)$,
where $\idn G$ is the independence number of $G$.
They also described some graphs for which equality
in~(1) is attained, and suggested that there might
be no other such graphs.

 Actually, (1) is one of a series of `linear
inequalities', which starts by
 $$\begin{array}{rcl}
 \ca E &\ge& 0, \\
 \ca E-\ca V+\idn G &\ge& 0, \\
 \ca E-3\ca V+5\idn G &\ge& 0, \hbox{ and} \\
 \ca E-5\ca V+10\idn G &\ge& 0 \,.
 \end{array}$$
 For each one of these `earlier' inequalities, the
graphs for which equality holds are edge number
critical (with respect to triangle freeness and to
vertex numbers and independence number upper
bounds), and they have been classified explicitly,
mainly by Radziszowski and Kreher; see e.~g.\
propositions~2.2 and~6.3 in \cite{RK}.

 \smallskip
 Some years later, I was able to confirm
Radziszowski's and Kreher's conjecture, and
included a proof in my book manuscript {\sl
Contributions to a Ramsey calculus\/} \cite B. Thus, the
inequality (1) is strict for all other triangle free
graphs. This consequence was quoted and employed in
2000 by Lesser \cite L.
 However, my still far from finished manuscript, and
thus the proof, have remained unpublished. It has
been pointed out to me that this is a non-optimal
state of matter. I therefore decided to present the
proof in a `stand-alone' article.

 However, just reproducing the proof from my
manuscript \cite B together with all its dependencies,
would amount to an unproportionally large article. I
am trying to make \cite B as self-contained as possible,
but on the other hand its proofs do contain numerous
internal references to more general results with
multiple applications.
 On the other hand, in the present article,
shortcuts are possible, largely due to the
possibility to refer to \cite{RK}, where in fact a
considerable part of the necessary ground work is
done.

 This article thus has a dual character. The
concepts, terminology, and notation largely follow
my manuscript, but the proofs as far as reasonably
possible are simplified by recycling the \cite{RK}
arguments and results. In particular, it in no ways
should depend on unpublished results, except in the
broader discussion in section~\ref{S:disc} at the
end, and in some of the footnotes in the earlier
sections (none of which contains facts used in the
proofs there).

 \bigskip
 \section{Fundamental concepts and notation.}
 \subsection{Basics.}
 Throughout this article, all considered graphs are
undirected, simple, and finite; thus, formally, a
graph is a pair $(V,E)$ of finite sets, such that
every element of $E$ is a 2-subset of $V$. In other
words, we demand that the cardinality $\ca V <
\infty$, and that
 $$E \subseteq {V \choose 2} = \{V \hbox{ subsets of
cardinality } 2\}.$$
 If $G = (V,E)$ we also let $V(G) = V$ and $E(G) =
E$.

 It will be convenient {\sl not\/} to disconsider
the empty graph, formally the pair $(\0,\0)$, but in
shorthand represented just by $\0$. Similarly, here
the natural numbers be $\N = \{0, 1, 2, 3,\ldots\}$
(including zero).

 When we consider a fixed graph $G = (V,E)$, the
shorthand notation often is extended to arbitrary
subsets of $V$; $W \subseteq V$ sometimes also may
be used for the induced subgraph $(W, E_W)$, where
$E_W$ in its turn is shorthand for $ E \cap {W
\choose 2}$. If $W$ and $X$ are two subsets of $V$,
then
 $$E_{W,X} = \bigl\{ \{w,x\} \in E : w \in W \land x
\in X \bigr\}\,.$$

 Recall that an {\sl isomorphism\/} between two
graphs $(V,E)$ and $(V',E')$ is a bijection between
$V$ and $V'$ which induces a bijection between $E$
and $E'$; let $G \simeq G'$ denote either an
isomorphism, or just the fact that $G$ and $G'$ are
isomorphic, depending on the context.
 If moreover \row vr and \row wr are sequences of
vertices in $G$ and $G'$, respectively, then
$(G,\row vr) \simeq (G',\row wr)$ denotes a {\sl
relative isomorphism}, i.~e., an isomorphism which
furthermore maps $v_i$ to $w_i$ for $ i =1 \kdots
r$, or the existence of such an isomorphism.

 There are numerous classes of graphs, which
technically are only defined up to isomorphisms. By
a slight abuse of terminology, often a graph will be
used, where more correctly an isomorphism class of
graphs should be treated. Likewise, $=$ (equal to)
may be used, where technically $\simeq$ (isomorphic
to) would be more correct.

 As usual, $K_i$ and \c i denote `the' complete
graph (or properly: a complete graph) and `the'
cycle graph on $i$ vertices, respectively. In
general, let $V(K_i) = \{\row ki\}$ and $V(\c i) =
\{\row ci\}$. \p i denotes `the' path graph with $i$
vertices (and thus $i-1$ edges). $K_{i,j}$ is the
complete bipartite graph with $i$ and $j$ vertices
in the respective parts. If different copies of
these or other graphs defined up to isomorphism are
needed in the same context, they are distinguished
by primes; as in \c5, $\c5'$, $\c5''$,~\dots\ for
different 5-cycles.

 In this article the {\sl sum\/} of two graphs is
their disjoint union; in other words, for $G =
(V,E)$ and $G' = (V',E')$ with $ V \cap V' = \0$, we
put $G+G' = (V \cup V', E \cup E')$; while, if $V
\cap V' \ne \0$, $G+G'$ only is defined up to
isomorphism, as the sum $G''+G'''$ for any $G''
\simeq G$ and $G''' \simeq G'$, such that $V(G'')
\cap V(G''') = \0$. For any $m \in \N$ and any graph
$G$, $mG$ is the sum of $m$ copies of $G$; in
particular, $0G = \0$ and $1G = G$.

 A graph $G = (V,E)$ is {\sl edge critical}, if the
removal of any edge increases the independence
number, i.e., if $\idn{(V,E')} > \idn G$ for every
proper subset $E'$ of $E$. If it is not edge
critical, then there is some edge $\ep \in E$, which
is {\sl redundant}, i.~e., such that $\idn
{(V,E\setminus\{\ep\})} = \idn G$.

 As usual, the {\sl distance\/} $\dist(v,w)$ between
two vertices is the smallest number of edges in any
path connecting them, if there is one, and $\infty$
else.

 The {\sl link\/} of a vertex $v$ in a graph $G =
(V,E)$ is (the induced subgraph on) the set of
vertices adjacent to $v$: $\lk v = \lk_G(v) =
\bigl\{w \in V : \{v,w\} \in E \bigr\} = \{w \in V :
\dist(v,w) = 1\}$\footnote{Graphs may be considered
as (the 1-skeletons of) flag abstract simplicial
complexes. From this point of view, the link of a
0-simplex $\{v\}$, in its usual sense, is precisely
$\lk v$.}. (Here and in the sequel, denoting the
graph may be omitted, if it is clear from the
context.) The (first) {\sl valency\/} of $v$ is $\v
v = \vind Gv = \ca {\lk v}$. The {\sl second
valency\/} of $v$ is the sum of all first valencies
of its neighbours:
 $$\vv2v = \vvind G2v = \v{\lk v} := \sum_{w\in\lk
v} \v w\,.$$

 For $v \in V$ and $d \in \N$, the $d$-{\sl
neighbourhood\/} of $v$ or the {\sl ball of radius
$d$ and centre $v$} is $\B dv = \B d{G,v} = \{w\in
V: \dist(v,w)\le d\}$. For $S = \{\row vr\}
\subseteq V$,
 $$\B dS = \B d{G,S} = \B d{\row vr} = \B d{G,\row
vr} = \bigcup\limits\3j1r \B d{v_j}\,.$$ A {\sl
monovalent\/} ({\sl bivalent, trivalent}, et cetera)
(vertex) is a vertex of valency 1 (2, 3, et cetera,
respectively).

 \subsection{Invariants.}
 Recall that a (proper) {\sl graph invariant\/} is a
number valued function $f$ on the set of all (finite
et cetera) graphs, such that
 $$G \simeq G' \implies f(G) = f(G') \,.$$
 The invariant $f$ is {\sl linear}, if in addition
it `respects sums', i.~e., if
 $$f(G+H) = f(G)+f(H), \ \forall\, G,H\,.$$
 In particular, then clearly $f(mG) = mf(G)$, and
$f(\0) = 0$.

 A {\sl generalised graph invariant\/} is defined
similarly, but some of the values may be
non-numbers. (We do not consider any kind of
linearity condition for non-proper invariants.)
An example is the {\sl girth},
 $$\girth G := \inf \{i : \o{\c i}G \ne 0\};$$
 thus, by the usual convention for infima of empty
sets of natural numbers, $\girth G = \infty$ if and
only if $G$ is acyclic, i.~e., is a forest.

 A {\sl linear inequality\/} is an inequality $f(G)
\ge g(G)$ involving two linear graph invariants $f$
and $g$, and which holds for all graphs of some
specified class, which is closed under addition.

 \smallskip
 For two graphs $H$ and $G$, let the {\sl number of
occurrences\/} of $H$ in $G$, $\o HG$, be the number
of $G$ subgraphs $G'$ (induced or not), such that
$G' \simeq H$. There is also a relative variant: For
$\row ur \in V(H)$ and $\row vr \in V(G)$, $\o
{H,\row ur} {G,\row vr}$ denotes the number of
subgraphs of $G$ which are isomorphic to $H$
relatively the $u_i$ mapping to the $v_i$. If there
is no possible ambiguity, the \row ur may be
omitted. Thus, $\o{\c i}{G,v}= \o{\c i,c_1}{G,v}$ or
$\o{\c i,c_1,c_2}{G,v,w}$ is the number of
$i$-cycles through a vertex $v$ or an edge $\{v,w\}$
in $G$, respectively.

 Note, that $\o HG$ is a graph invariant (with
respect to $G$), for any fixed $H$. This invariant
is linear if and only if $H$ is connected. Two
such linear invariants are
 $$n(G) = \o {K_1} G = \ca {V(G)} \ \hbox{ and }\
 e(G) = \o {K_2}G = \ca {E(G)} \,.$$

 For each natural number $d$, the number of
$d$-valents is a graph invariant, denoted \nv dG. 

 Another important linear graph invariant is the
{\sl independence number}, the maximal size of an
independent subset of vertices:
 $$\idn G = \max \; \bigl( \ca S : S \subseteq V
\land {S \choose 2} \cap E = \0 \bigr) \,.$$
 Likewise, the number $\comp G$ of (connected)
components of $G$ is a linear graph invariant.
 On the other hand, the {\sl clique number\/}
 $\cqn G = \max \,(i : \o {K_i} G > 0)$ is a graph
invariant, but not linear. In fact, $\cqn{G+H} =
\max\, \bigl( \cqn G, \cqn H \bigr)$. ($\cqn \0 =
0$.)

 The graph $G$ is {\sl triangle free}, if $\cqn G
\le 2$, or, equivalently, if $\girth G \ge 4$. In
the later sections of this article, we only consider
triangle free graphs. The triangle free graph $G$ is
{\sl square free}, if in addition it does not
contain any 4-cycle (or ``square''), or,
equivalently, if $\girth G \ge 5$.

 Directly from the definitions we get
 \begin{lemma}
 \label{L:21}
A linear combination of linear
graph invariants is itself a linear graph
invariant.\noproof
 \end{lemma}

 In this article, the two most important linear
invariants formed as linear combinations are
 $$t(G) := \ell_6(G) = e(G) - 6n(G) + 13 \idn G
\quad \hbox{ and }\quad
 q(G) := \ell_5(G) = e(G) - 5n(G) + 10\idn G \,.$$

 \subsection{Independence stability.}
 A {\sl destabilising subset\/} or {\sl
destabiliser\/} $M$ in a graph $G = (V,E)$ is a
subset of $V$, such that the induced subgraph on $V
\setminus M$ has a lower independence number than
$\idn G$. The graph $G$ is {\sl $s$-stable}, if it
has no destabiliser of size $\le s$, and is {\sl
strongly\/} $s$-stable, if in addition no
destabilising subset of size $s+1$ is independent.
(Often, the induced subgraph on $M$ also is called
$M$.) Now, if $v \in V = V(G)$, and $S$ is an
independent subset of $V \setminus \B1v$, then $S
\cup \{v\}$ also is independent; whence we directly
get
 \begin{lemma}
 \label{L:22}
 For any vertex $v$ in a graph
$G$, $B(G,v;1)$ destabilises $G$.\noproof
 \end{lemma}

 If $S = \{\row vr\}$ is an independent subset of $G
= (V,E)$, then
 $$G_{\row vr} = (V_{\row vr}, E_{\row vr})$$
 denotes the induced subgraph on everything but the
neighbourhood of $S$; in other words, $V_{\row vr} =
V \setminus \B1{G,S}$ and $E_{\row vr} = E \cap
{V_{\row vr} \choose 2}$. (Whenever we employ the
notation $G_{\row vr}$, we indeed assume that the
$v_i$ are different and form an independent set.)
Since obviously
 $$G_{\row vr} =
(\cdots((G_{v_1})_{v_2})\cdots)_{v_r},$$
 and by inductive use of lemma~2.2, we have
 \begin{lemma}
 \label{L:23}
 $\idn {G_{\row vr}} \le \idn G-r$.\noproof
 \end{lemma}

 In some interesting situations, we have equality.
 \begin{lemma}
 \label{L:ec1um}
 If $G = (V,E)$ is edge critical and $v \in V$, then
$\idn{G_v} = \idn G-1$.
 \end{lemma}
 \begin{proof}
 This is immediate from the linearity, if $\v v=0$,
since then $G = \{v\} + (V\setminus\{v\},E)$.

 Else, choose a $w \in \lk v$, and let $k = \idn G$
and $E' = E \setminus \bigl\{ \{v,w\} \bigr\}
\subset E$. By the edge criticality, there is a
$(k+1)$-subset $S$ of $V$, which is independent in
$(V,E')$. Now, if $S$ did not contain both $v$ and
$w$, it were independent in $G$ as well, against the
assumptions. Thus, instead, $v,w \in S$, and $S
\setminus \{v,w\}$ is an independent $(k-1)$-subset
of $G_v$.

 Thus, $\idn{G_v} \ge k-1 = \idn G-1$, whence we
indeed have equality by lemma~\ref{L:23}.
 \end{proof}

 Obviously, $n(G_v) = n(G) - \v v-1$. Moreover, if
$\o{K_3,k_1}{G,v} = 0$, then $e(G_v) = e(G) - \vv2v$
(since then \vv2v counts each edge in $E \setminus
E_v$ exactly once). Generalising, we get
 \begin{lemma}
 \label{L:stsg}
 If $G = (V,E)$ is a triangle free graph, and $S =
\{\row vr\} \subseteq V$ an independent set, then
 $$n(G) - n(G_{\row vr}) = \ca{\B1S} = r + \sum\3i1r
\vind {G_{\row v{i-1}}}{v_i}\,,\ \hbox{ and}$$
 $$e(G) - e(G_{\row vr}) = \v{\B1S} - \ca E_{\B1S} =
\sum\3i1r \vvind {G_{\row v{i-1}}}2{v_i}\,.
\eqno\square$$
 \end{lemma}

 \begin{lemma}
 \label{L:prpt}
 Let $G$ be an edge critical and connected triangle
free graph, and let $v \in V(G)$.

 $(a)$ If $v$ is bivalent, then $\comp {G_v} = 1$.

 $(b)$ If $v$ is trivalent, $G_v = G'+G''$, $G' \neq
\0 \ne G''$, $G'$ is strongly $s'$-stable, and $G''$
is strongly $s''$-stable, then
 $$\vv2v \ge s'+s''+6 \,.$$
 \end{lemma}

 \begin{proof}
 First, note that
 $$\idn {G_v} = \idn G-1 \ge \v v-1 \ge 1 \implies
G_v \ne \0,$$
 by lemma~\ref{L:ec1um}.
 Thus, in case~$(a)$, for a contradiction, we also
may assume $G_v = G'+G''$ with both $G^{(\nu)}$
non-empty.

 In both cases, let $\lk v = \{\row w{\v v}\}$, $X_i
= \lk(w_i) \setminus \{v\}$ (for $i = 1 \kdots \v
v$), $X = \bigcup_{i=1}^{\v v} X_i$, $X^{(\nu)} = X
\cap V(G^{(\nu)}) = \B2v \setminus \B1v$ (for $\nu =
1,2$), and $X_i^{(\nu)} = X_i \cap V(G^{(\nu)})$. By
connectedness, both $E_{\lk(v),X^{(\nu)}}$ must be
non-empty.

 Now, in case $(a)$, if $X'$ did not destabilise
$G'$, then the edges in $E_{\lk(v),X'}$ were
redundant, and else so were the edges in
$E_{\lk(v),X''}$, in either case contradicting the
edge criticality.

 In case~$(b)$, similarly, $X^{(\nu)}$ must
destabilise $G^{(\nu)}$ for both $\nu$. By the
assumptions, in particular,
 $$a^{(\nu)} := \ca{E_{X^{(\nu)},\lk v}} =
\ca{X_1^{(\nu)}} + \ca{X_2^{(\nu)}} +
\ca{X_3^{(\nu)}} \ge \ca{X^{(\nu)}} \ge s^{(\nu)}+1,
\quad \nu=1,2\,. \leqno (2)$$
 In particular, clearly $\vv2v = 3 + \ca{E_{X,\lk
v}} = 3+a'+a'' \ge s'+s''+5$, and it suffices to
prove that equality would lead to a contradiction.

 Indeed, the only way to have equality would be to
have equalities in~(2), too, whence on the one hand
all the $X_i^{(\nu)}$ were disjoint, while on the
other hand $X^{(\nu)}$ were a minimal and
non-independent destabiliser of $G^{(\nu)}$, for
both $\nu$. In particular, for each $\nu$, at least
two of $X_1^{(\nu)}$, $X_2^{(\nu)}$, and
$X_3^{(\nu)}$ were non-empty; whence, without loss
of generality, we could assume $X'_3 \ne \0 \ne
X_3''$. However, then the $X_1^{(\nu)} \cup
X_2^{(\nu)}$ were proper subsets of the $X^{(\nu)}$,
whence neither $X_1^{(\nu)} \cup X_2^{(\nu)}$ would
destabilise $G_v$, whence nor would $X_1 \cup X_2$.
Thus, there were an independent $(k-1)$-subset $S$
of $V(G_v) \setminus (X_1 \cup X_2)$; but then $S
\cup \{w_1,w_2\}$ were an independent $(k+1)$-subset
of $G$, against the assumptions; the sought
contradiction.
 \end{proof}

 \subsection{E-numbers.}
 An $(i,j;n,e)$ {\sl realiser\/} is a graph $G$ with
$\cqn G < i$, $\idn G < j$, $n(G) = n$, and $e(G) =
e$. The {\sl e-number\/} \e i,j;n \ is the minimal
number $e$, such that there are $(i,j;n,e)$
realisers, or $\infty$, if no such realisers exist
for any $e$. These numbers are closely related to
the (classical 2-colours) Ramsey numbers; in fact,
these may be defined by
 $$\rr ij := \min \,\{ n \in \N : \e i,j;n = \infty
\}\,.$$
 There has been some efforts to determine the
e-numbers \e3,k+1;n . In \cite{GR} {\sl inter
alia\/} all such e-numbers for $k+1 \le 10$ are
listed. The non-negative linear invariants and
constructions of graphs with invariant value 0
provide some infinite families of e-values; as shown
in \cite[theorems~1 and~4]{RK'} and
\cite[theorem~5.1.1 and corollary~5.3.4]{RK},
 \begin{proposition}[Radziszowski, Kreher]
 \label{P:eval}
 For $k \ge 4$ and either $0 \le n \le 3.25k-1$ or
$n = 3.25k$,
 $$\e3,k+1;n = \max\, (0, n-k, 3n-5k, 5n-10,
6n-13k)\,.$$
 Moreover, for all $n$ and $k$, $\e3,k+1;n \ge
6n-13k$.
 \noproof
 \end{proposition}

 \section{Constructing graphs step by step.}
 \label{S:step}
 From now on, all graphs considered in this article
are triangle free, unless explicitly otherwise
denoted.

 Graphs `close to a limit' will tend to share more
structure than `graphs in general'. Often, they also
have subgraphs `close to' that limit. This may make
their structure tractible to recursive treatment.

 Radziszowski and Kreher consider {\sl
$d$-extensions\/} $G_v \subset G$, where $\v v=d$,
and where $G_v$ is demanded to be edge number
critical with respect to $n(G_v)$ and $\idn{G_v}$.
Here, we consider somewhat more general kinds of
extensions, called {\sl stitches}. Indeed, stress is
both on the details for each such stitch (or
extension), and on the way the whole graph may be
composed by such steps. I found the analogy with
knitting or chrochet rather apt, and therefore
partly follow chochet terminology.

 Formally, quite generally, a stitch is a pair
$(G,G') = \bigl( (V,E),\;\allowbreak (V',E') \bigr)$
of graphs, such that $G'$ is the induced subgraph of
$G$ on $V' \subset V$ and $M = V \setminus V'$ is a
minimal destabiliser of $G$. The latter condition
precisely means that
 $$\idn {G'}+1 = \idn {G''} = \idn G $$
 for the induced graph $G''$ on any $V''$ with $V'
\subset V'' \subseteq V$. Less formally, this $G$
may be called `a stitch of $G'$'. The stitch is
classified by the corresponding minimal destabiliser
$M$, and by the way $M$ is `fastened' at $G'$.

 Thus, generalising the concept $d$-extensions from
\cite{RK}, a {\sl $d$-stitch\/} is a graph $G$ and a
specified vertex $v \in V(G)$, such that $\v v=d$
and $\idn{G_v} = \idn G-1$, but that no edge in \B1v
may be removed without increasing the independence
number. Note, that the last condition certainly is
satisfied, if $G$ is edge critical.

The most simple case is the 1-stitch, extending $\0$
to \p2, the 2-vertices path.

 In general, the graphs we are interested of here
may be constructed from scratch by a number of
$d$-stitches, for various $d$.

 If $H$ is a $d$-stitch of $H_v$, then call
$v$ the {\sl apex\/} of the stitch. Moreover, the
set $X$ of vertices of distance exactly 2 from $v$
is called the {\sl base} of the stitch, which
also is said to be {\sl based} at $X$. Note, that
$X$ is a subset of $V(H_v)$. If moreover the
stitch is a 2-stitch, then $X$ destabilises
$H_v$, and is bipartite (since the neighbourhoods of
the two neighbours of $v$ must be an independent
set).

 Conversely, a bipartite minimal destabiliser $M$ of
a graph $G$ can be used as base of a 2-stitch of
$G$. (Note, that the `bipartitivity' is
self-evident, if $M$ contains at most four vertices,
since that is too few for a 5-cycle or a larger
cycle of odd length.) If $G$ is edge critical, and
$v \in V$, then \B1v is such a destabiliser. In this
case, the corresponding 2-stitch also is said to be
{\sl based at $v$}. In these cases, the 2-stitch
(which is uniquely defined up to isomorphisms) also
may be denoted $\crt_2(G;M)$ or $\crt_2(G;v)$; the
index 2 may be omitted. The apex of that stitch, say
$v'$, may be used as the base of a new stitch
$\crt\bigl(\crt(G;M);v'\bigr)$, also denoted
$\crt^2(G;M)$; {\sl et cetera}.

 We now may form the most fundamental crochet, a
(simple) {\sl chain}, by successively adding
2-stitches to \c5\footnote{Actually, $\c5 \simeq
\crt_2 (\p2)$, and $\p2 \simeq \crt_1 (\0)$, whence
instead indeed we might start ``from scratch'',
putting $\w1 := \p2$ and $\w0 := \0$; however, we
have no use of this notation in this article.}.
Thus, put $\w2 := \c5$, and $\w k := \crt^{k-2}
(\c5;c_1) \simeq \crt(\w{k-1};p)$ for any $k \ge 3$
(and any vertex $c_1$ in \c5 or bivalent $p$ in
\w{k-1}, respectively).

 \bigbreak
 \begin{picture}(370,130)

 \thicklines

  \put(10,110){\w3, alias $\c5\,$.}

 \put(25,50){\circle*5}
 \put(16,77){\circle*5}
 \put(40,95){\circle*5}
 \put(64,77){\circle*5}
 \put(55,50){\circle*5}

 \put(25,50){\line(-1,3)9}
 \put(16,77){\line(4,3){24}}
 \put(40,95){\line(4,-3){24}}
 \put(64,77){\line(-1,-3)9}
 \put(55,50){\line(-1,0){30}}

 \thinlines
 \put(-10,0){Its crochet pattern.}

 \put(30,20){\circle*3}
 \put(50,20){\circle*3}
 \put(30,20){\line(1,0){20}}

 \thicklines
  \put(200,110){\w5 as a 2-stitch of $\w4\,$.}

 \put(200,80){\circle*5}
 \put(230,80){\circle*5}
 \put(260,80){\circle*5}
 \put(290,80){\circle*5}
 \put(320,80){\circle*5}
 \put(350,80){\circle*5}
 \put(200,50){\circle*5}
 \put(230,50){\circle*5}
 \put(275,50){\circle*5}
 \put(320,50){\circle*5}
 \put(350,50){\circle*5}

 \put(201,80){\line(1,0){148}}
 \put(201,50){\line(1,0){148}}
 \put(200,50){\line(0,1){30}}
 \put(230,50){\line(1,1){30}}
 \put(275,50){\line(-3,2){45}}
 \put(275,50){\line(3,2){45}}
 \put(320,50){\line(-1,1){30}}
 \put(350,50){\line(0,1){30}}

 \thinlines
 \put(140,80){\circle5}
 \put(170,80){\circle5}
 \put(170,50){\circle5}

 \put(142,80){\line(1,0){26}}
 \put(172,80){\line(1,0){26}}
 \put(172,50){\line(1,0){26}}
 \put(142,78){\line(1,-1){26}}
 \put(172,79){\line(2,-1){57}}

 \thinlines
 \put(220,0){Its crochet pattern.}

 \put(225,20){\circle3}
 \put(245,20){\circle*3}
 \put(265,20){\circle*3}
 \put(285,20){\circle*3}
 \put(305,20){\circle*3}

 \put(245,20){\line(1,0){60}}

 \put(226.3,20){\line(1,0){18}}

 \end{picture}
 \bigbreak

 Next, we may consider 3-stitches of chains. The
simplest such are the {\sl bicycles}. The bicycle
\DW k ($k \ge 4$) consists of an (induced) {\sl
outer cycle\/} $\{\row d{2k}\} \simeq \c{2k}$ and an
(induced) {\sl inner cycle\/} $\{\row ek\} \simeq \c
k$, with the connecting edges $\{d_{2i-2,e_i}\}$ and
$\{d_{2i+1},e_i\}$ for $i=1\kdots k$ (where as usual
the outer and inner cycle vertex indices may be
calculated modulo $2k$ and modulo $k$,
respectively).

 \begin{figure*}[t]
 \begin{picture}(350,170)

 \thicklines
 \put(35,5){The bicycle $\DW6\,$.}

 \put(45,95){\circle*5}
 \put(60,120){\circle*5}
 \put(90,120){\circle*5}
 \put(105,95){\circle*5}
 \put(90,70){\circle*5}
 \put(60,70){\circle*5}

 \put(20,110){\circle*5}
 \put(35,135){\circle*5}
 \put(60,150){\circle*5}
 \put(90,150){\circle*5}
 \put(115,135){\circle*5}
 \put(130,110){\circle*5}
 \put(130,80){\circle*5}
 \put(115,55){\circle*5}
 \put(90,40){\circle*5}
 \put(60,40){\circle*5}
 \put(35,55){\circle*5}
 \put(20,80){\circle*5}

 \put(45,95){\line(3,5){15}}
 \put(60,120){\line(1,0){30}}
 \put(90,120){\line(3,-5){15}}
 \put(105,95){\line(-3,-5){15}}
 \put(90,70){\line(-1,0){30}}
 \put(60,70){\line(-3,5){15}}

 \put(20,110){\line(3,5){15}}
 \put(35,135){\line(5,3){25}}
 \put(60,150){\line(1,0){30}}
 \put(90,150){\line(5,-3){25}}
 \put(115,135){\line(3,-5){15}}
 \put(130,110){\line(0,-1){30}}
 \put(130,80){\line(-3,-5){15}}
 \put(115,55){\line(-5,-3){25}}
 \put(90,40){\line(-1,0){30}}
 \put(60,40){\line(-5,3){25}}
 \put(35,55){\line(-3,5){15}}
 \put(20,80){\line(0,1){30}}

 \put(45,95){\line(-1,4){10}}
 \put(60,120){\line(1,1){30}}
 \put(90,120){\line(4,-1){40}}
 \put(105,95){\line(1,-4){10}}
 \put(90,70){\line(-1,-1){30}}
 \put(60,70){\line(-4,1){40}}
 \put(45,95){\line(-1,-4){10}}
 \put(60,120){\line(-4,-1){40}}
 \put(90,120){\line(-1,1){30}}
 \put(105,95){\line(1,4){10}}
 \put(90,70){\line(4,1){40}}
 \put(60,70){\line(1,-1){30}}

 \put(49,93){$e_1$}
 \put(60,110){$e_2$}
 \put(80,110){$e_3$}
 \put(90,93){$e_4$}
 \put(80,77){$e_5$}
 \put(60,77){$e_6$}

 \put(05,78){$d_1$}
 \put(05,108){$d_2$}
 \put(21,137){$d_3$}
 \put(55,157){$d_4$}
 \put(86,157){$d_5$}
 \put(118,137){$d_6$}
 \put(134,108){$d_7$}
 \put(134,78){$d_8$}
 \put(118,49){$d_9$}
 \put(85,27){$d_{10}$}
 \put(54,27){$d_{11}$}
 \put(20,48){$d_{12}$}

 \thinlines
 \put(200,5){Its crochet pattern.}

 \put(235,104){\circle*3}
 \put(250,113){\circle*3}
 \put(265,104){\circle*3}
 \put(265,86){\circle*3}
 \put(250,77){\circle*3}
 \put(235,86){\circle*3}

 \put(235,104){\line(5,3){15}}
 \put(250,113){\line(5,-3){15}}
 \put(265,104){\line(0,-1){18}}
 \put(265,86){\line(-5,-3){15}}
 \put(250,77){\line(-5,3){15}}
 \put(235,86){\line(0,1){18}}

 \end{picture}
 \end{figure*}

 Indeed, for any trivalent $v$ in \DW k, $(\DW k)_v
\simeq \w{k-1}$.

 \medbreak
 The bicycles may be considered as crochet loops.
They are rather symmetric, and like the chains
well-defined up to isomorphism, by means of the
single parameter $k = \idn{\DW k}$. However, the
next class we consider are `loop-chains', where we
add a succession of 2-stitches to a bicycle, to get
a pending attached chain; and here we need three
parameters: The length of the loop, the length of
the chain, and a description of how the chain is
attached to the loop.\footnote{In the crochet
pattern, the attachment is represented by a
trivalent, while the mood of attachment conveniently
may be represented by associating outgoing
directions to some of the edges at that trivalent.}
Here, we only are interested in attachments at
destabilisers of minimal size; i.~e., a loop-chain
may be written as $\crt^{k-l}(\DW l;M)$, where $l
\ge 4$, $k-l \ge 1$, and $M$ is a \DW l destabiliser
of size four. As we shall see in lemma~\ref{L:B},
there are up to isomorphism just three possible $M$,
isomorphic to $K_{1,3}$, \p4, or \c4, respectively,
whence we may let the third parameter range over
just the set $\{K_{1,3}, \p4, \c4\}$.

 \bigbreak
 \begin{picture}(370,155)



 \thicklines
 \put(-5,135){The loop-chain $\crt^2(\DW4;\p4)$.}

 \put(30,100){\circle*5}
 \put(50,100){\circle*5}
 \put(50,80){\circle*5}
 \put(30,80){\circle*5}

 \put(10,100){\circle*5}
 \put(30,120){\circle*5}
 \put(50,120){\circle*5}
 \put(70,100){\circle*5}
 \put(70,80){\circle*5}
 \put(50,60){\circle*5}
 \put(30,60){\circle*5}
 \put(10,80){\circle*5}

 \put(30,100){\line(1,0){20}}
 \put(50,100){\line(0,-1){20}}
 \put(50,80){\line(-1,0){20}}
 \put(30,80){\line(0,1){20}}

 \put(10,100){\line(1,1){20}}
 \put(30,120){\line(1,0){20}}
 \put(50,120){\line(1,-1){20}}
 \put(70,100){\line(0,-1){20}}
 \put(70,80){\line(-1,-1){20}}
 \put(50,60){\line(-1,0){20}}
 \put(30,60){\line(-1,1){20}}
 \put(10,80){\line(0,1){20}}

 \put(10,100){\line(1,-1){40}}
 \put(50,120){\line(-1,-1){40}}
 \put(70,80){\line(-1,1){40}}
 \put(30,60){\line(1,1){40}}

 \put(90,100){\circle*5}
 \put(110,100){\circle*5}
 \put(130,100){\circle*5}
 \put(90,80){\circle*5}
 \put(110,80){\circle*5}
 \put(130,80){\circle*5}

 \put(90,100){\line(-2,1){40}}
 \put(90,100){\line(-1,-1){20}}
 \put(90,100){\line(1,-1){20}}
 \put(90,100){\line(1,0){40}}
 \put(90,80){\line(-1,1){20}}
 \put(90,80){\line(-2,-1){40}}
 \put(90,80){\line(1,1){20}}
 \put(90,80){\line(1,0){40}}
 \put(130,80){\line(0,1){20}}

 \thinlines
 \put(10,0){Its crochet pattern.}

 \put(44,20){\circle*3}
 \put(32,32){\circle*3}
 \put(44,44){\circle*3}
 \put(56,32){\circle*3}
 \put(73,32){\circle*3}
 \put(90,32){\circle*3}

 \put(32,32){\line(1,1){12}}
 \put(32,32){\line(1,-1){12}}
 \put(56,32){\vector(-1,1){11}}
 \put(56,32){\vector(-1,-1){11}}
 \put(56,32){\vector(1,0){15.5}}
 \put(73,32){\line(1,0){17}}


 \thicklines
  \put(180,135){The loop-chain
$\crt^3(\DW4;K_{1,3})$.}

 \put(200,100){\circle*5}
 \put(220,100){\circle*5}
 \put(220,80){\circle*5}
 \put(200,80){\circle*5}

 \put(180,100){\circle*5}
 \put(200,120){\circle*5}
 \put(220,120){\circle*5}
 \put(240,100){\circle*5}
 \put(240,80){\circle*5}
 \put(220,60){\circle*5}
 \put(200,60){\circle*5}
 \put(180,80){\circle*5}

 \put(200,100){\line(1,0){20}}
 \put(220,100){\line(0,-1){20}}
 \put(220,80){\line(-1,0){20}}
 \put(200,80){\line(0,1){20}}

 \put(180,100){\line(1,1){20}}
 \put(200,120){\line(1,0){20}}
 \put(220,120){\line(1,-1){20}}
 \put(240,100){\line(0,-1){20}}
 \put(240,80){\line(-1,-1){20}}
 \put(220,60){\line(-1,0){20}}
 \put(200,60){\line(-1,1){20}}
 \put(180,80){\line(0,1){20}}

 \put(180,100){\line(1,-1){40}}
 \put(220,120){\line(-1,-1){40}}
 \put(240,80){\line(-1,1){40}}
 \put(200,60){\line(1,1){40}}

 \put(260,100){\circle*5}
 \put(280,100){\circle*5}
 \put(300,100){\circle*5}
 \put(320,100){\circle*5}
 \put(340,100){\circle*5}
 \put(260,80){\circle*5}
 \put(290,80){\circle*5}
 \put(320,80){\circle*5}
 \put(340,80){\circle*5}

 \put(260,80){\line(-1,1){20}}
 \put(260,80){\line(-2,1){40}}
 \put(260,80){\line(-2,-1){40}}
 \put(260,80){\line(1,0){80}}
 \put(260,80){\line(1,1){20}}
 \put(290,80){\line(-3,2){30}}
 \put(290,80){\line(3,2){30}}
 \put(320,80){\line(-1,1){20}}
 \put(340,80){\line(0,1){20}}
 \put(260,100){\line(-1,-1){20}}
 \put(260,100){\line(1,0){80}}

 \thinlines
 \put(215,0){Its crochet pattern.}

 \put(234,20){\circle*3}
 \put(222,32){\circle*3}
 \put(234,44){\circle*3}
 \put(246,32){\circle*3}
 \put(263,32){\circle*3}
 \put(280,32){\circle*3}
 \put(297,32){\circle*3}

 \put(222,32){\line(1,1){12}}
 \put(222,32){\line(1,-1){12}}
 \put(246,32){\line(-1,1){11}}
 \put(246,32){\vector(-1,-1){11}}
 \put(246,32){\line(1,0){51}}

 \end{picture}
 \bigbreak

 \bigbreak
 \begin{picture}(350,210)

 \thicklines
 \put(0,5){The loop-chain $\crt(\DW6;C_4)$.}

 \put(45,95){\circle*5}
 \put(60,120){\circle*5}
 \put(90,120){\circle*5}
 \put(105,95){\circle*5}
 \put(90,70){\circle*5}
 \put(60,70){\circle*5}

 \put(20,110){\circle*5}
 \put(35,135){\circle*5}
 \put(60,150){\circle*5}
 \put(90,150){\circle*5}
 \put(115,135){\circle*5}
 \put(130,110){\circle*5}
 \put(130,80){\circle*5}
 \put(115,55){\circle*5}
 \put(90,40){\circle*5}
 \put(60,40){\circle*5}
 \put(35,55){\circle*5}
 \put(20,80){\circle*5}

 \put(45,95){\line(3,5){15}}
 \put(60,120){\line(1,0){30}}
 \put(90,120){\line(3,-5){15}}
 \put(105,95){\line(-3,-5){15}}
 \put(90,70){\line(-1,0){30}}
 \put(60,70){\line(-3,5){15}}

 \put(20,110){\line(3,5){15}}
 \put(35,135){\line(5,3){25}}
 \put(60,150){\line(1,0){30}}
 \put(90,150){\line(5,-3){25}}
 \put(115,135){\line(3,-5){15}}
 \put(130,110){\line(0,-1){30}}
 \put(130,80){\line(-3,-5){15}}
 \put(115,55){\line(-5,-3){25}}
 \put(90,40){\line(-1,0){30}}
 \put(60,40){\line(-5,3){25}}
 \put(35,55){\line(-3,5){15}}
 \put(20,80){\line(0,1){30}}

 \put(45,95){\line(-1,4){10}}
 \put(60,120){\line(1,1){30}}
 \put(90,120){\line(4,-1){40}}
 \put(105,95){\line(1,-4){10}}
 \put(90,70){\line(-1,-1){30}}
 \put(60,70){\line(-4,1){40}}
 \put(45,95){\line(-1,-4){10}}
 \put(60,120){\line(-4,-1){40}}
 \put(90,120){\line(-1,1){30}}
 \put(105,95){\line(1,4){10}}
 \put(90,70){\line(4,1){40}}
 \put(60,70){\line(1,-1){30}}

 \put(60,180){\circle*5}
 \put(90,180){\circle*5}
 \put(75,195){\circle*5}
 \put(75,195){\line(1,-1){15}}
 \put(75,195){\line(-1,-1){15}}
 \put(60,180){\line(1,-1){30}}
 \put(60,180){\line(1,-2){30}}
 \put(90,180){\line(-1,-1){30}}
 \put(90,180){\line(-1,-2){30}}

 \thinlines
 \put(200,5){Its crochet pattern.}

 \put(235,104){\circle*3}
 \put(250,113){\circle*3}
 \put(265,104){\circle*3}
 \put(265,86){\circle*3}
 \put(250,77){\circle*3}
 \put(235,86){\circle*3}
 \put(250,131){\circle*3}

 \put(235,104){\line(5,3){15}}
 \put(250,113){\line(5,-3){15}}
 \put(265,104){\line(0,-1){18}}
 \put(265,86){\line(-5,-3){15}}
 \put(250,77){\line(-5,3){15}}
 \put(235,86){\line(0,1){18}}
 \put(250,113){\vector(0,1){16.5}}

 \end{picture}
 \bigbreak

 Finally, we consider one famous ``cyclic graph'',
the unique triangle free graph on 13 vertices with
independence number 4. The graph often is called
$H_{13}$; here, I use the more descriptive notation
\5, signifying that it has 13 indexed vertices, and
that each vertex is connected to those with a
difference 1 or 5 in indices (again calculated
modularly, where appropriate)\footnote{Another,
more algebraic way to define \5 is by taking the
Galois field $GF(13)$ as its vertex set, and putting
$E(\5) := \bigl\{ \{x,y\} \subset GF(13) : (x-y)^4 =
1 \bigr\}$.}. By the definitions and inspection, \5
is a 4-stitch of \w3, where any vertex in \5 may be
chosen as the apex.

 \bigskip
\section{The main result.}
\label{thm}
 Recall that $q(G) = e(G)-5n(G)+10\idn G$ and $t(G)
= e(G)-6n(G)+13\idn G$. The following is a summary
of results by Radziszowski and Kreher, collected
from \cite{RK}:
 \begin{proposition}[Radziszowski, Kreher]
 \label{P:41}
 For any triangle free graph $G$, $q(G) \ge 0$ and
$t(G) \ge 0$. Moreover, $q(G) = 0$ if and only if
$G$ is a sum of chains and bicycles. If $G$ is a sum
of bicycles and copies of \5, then $t(G) = 0$; and
as a partial converse, if $\idn G \le 6$ and
$t(G)=0$, then indeed $G$ is such a sum. Finally, in
any case, if $t(G)=0$ but $G \ne \0$, then $3 \le
\mval G $, but $G$ is not 3-regular.
 \end{proposition}
 The $q$ inequality and extremal graph
characterisation is proposition~2.2~(b) and
theorem~4.3.1 in \cite{RK}, respectively. The $t$
inequality is their (main) theorem~5.1.1, the
equality cases is proposition~6.3~(e), and the $\idn
G \le 6$ graphs characterisation is a remark at the
beginning of the proof of lemma~5.1.5, in p.~77.
 That $\mval G \ge 3$ is contained in lemma~5.1.5.
 The last statement, that $G$ cannot be non-empty
but 3-regular, is implicitly noted {\sl en
passant\/} in the proof of lemma~5.1.6, since such a
$G$ would be a ``minimum graph with average degree
not exceeding 10/3'', in the terminology of that
proof.

 They suggest that $t(G)=0$ {\sl only\/} for sums
of bicycles and \5 copies. This indeed is true, but
for my proof to work, I also had to characterise
some triangle free graphs with $t(G)=1$:
 \begin{theorem}
 \label{T:1}
 Let $G$ be any triangle free graph.
 \clause{$(a)$}{If $t(G)=0$, then each component of
$G$ is a bicycle or isomorphic to \5.}
 \clause{$(b)$}{If $t(G) = 1$ and $\mval G \le 2$,
then one component of $G$ is a chain or a \loch, and
any other components are bicycles and copies of \5.}
 \end{theorem}

Put
 $$\begin{array}{rcl}
 \Psi &=& \{G : \cqn G < 3 \land t(G)=0\},\\
 \Psi' &=& \{G : \cqn G < 3 \land t(G)=1 \land \mval
G \le 2\},\\
 \Gamma &=& \{G : \hbox{ each component of $G$ is a
bicycle or a \5 copy$\}$, and}\\
 \Gamma' &=& (\{ \hbox{chains of lengths } \ge2\}
\cup \{\hbox{loop-chains}\}) + \Gamma,
 \end{array}$$
 with the usual algebraic interpretation of a sum of
sets of addable elements. The theorem now may be
reformulated as
 $$\Psi = \Gamma \hbox{ and } \Psi' = \Gamma'.
\leqno(3)$$
 By just calculating $t(G)$ and $\mval G $ for the
connected members of $\Psi \cup \Psi'$ and employing
linearity, it is easy to see that indeed $\Gamma
\subseteq \Psi$ (as noted in
proposition~\ref{P:41}), and that $\Gamma' \subseteq
\Psi'$. Thus, we only have to prove the converse
inclusions in~(3).

 \smallskip
 Since $t$ is linear, and is non-negative on
triangle free graphs, the class $\Psi$ is closed
under addition; in fact, a graph $G$ belongs to
$\Psi$ if and only if every component of $G$ does.
On the other hand, $\Psi'$ obviously is not closed
under addition. Instead, every graph in $\Psi'$ has
exactly one component $C$ with $t(C)=1$, while the
other components belong to $\Psi$. Moreover, since
every one of the latter components has minimal
valency $\ge3$ by proposition~\ref{P:41}, but $G$
does not, in fact $C \in \Psi'$. In other words,
 $$\Psi' = \{G \in \Psi' : \comp G = 1\} + \Psi.
\leqno(4)$$
 Thus, while $\Psi \cup \Psi'$ is not closed under
addition, it is closed under taking components or
other summands:
 $$G'+G'' \in \Psi \cup \Psi' \implies G',G'' \in
\Psi \cup \Psi'\,. \leqno(5)$$
 Thus, in order to prove the theorem, it is
sufficient to prove that a connected $G$ in $\Psi$
($\Psi'$) also must belong to $\Gamma$ ($\Gamma'$,
respectively), by means of induction with respect to
$\idn G$. This will be done in section~\ref{S:prf}.

 \smallskip
 As a direct consequence of the theorem, and since
there are realisers for $(3,6;16,32)$, $(3,7;19,37)$
and $(3,8;22,42)$ and by linearity, we get a slight
improvement of proposition~\ref{P:eval}:
 \begin{corollary}
 \label{cor}
 Let $n$ and $k \ge 5$ be
integers. If $3.25k-1 < n < 3.25k$, then $\e3,k+1;n
= 6n-13k+1$, and if $n > 3.25k$, then $\e3,k+1;n \ge
6n-13k+1$.\noproof
 \end{corollary}
 
 \bigskip \bigskip
 \section{Preparatory results.}
 I'll start by collecting the further needed results
in a few lemmata. As far as possible, I refer their
proofs to corresponding \cite{RK} results.

 \bigskip
 \begin{lemma}
 \label{L:A}
 For any $k \ge 2$, $\idn{\w k}=k$, and \w k is
2-stable and has no destabilising subset of size 3
other than \B1v for any bivalent $v$ therein.
 \end{lemma}

 \begin{proof}
 This is essentially a reformulation of
\cite[lemma~4.2.2~(a)]{RK}.
 \end{proof}

 \bigskip
 \begin{lemma}
 \label{L:B}
 For any $k \ge 4$, $\idn{\DW k}=k$, and \DW k is
3-stable and only has three kinds of destabilising
4-sets of vertices, namely
 \hfill\break
 (1): A ball \B1{d_i}, i.~e., a trivalent, together
with its three neighbours;
 \hfill\break
 (2): an induced path $\{d_{2d-1}, d_{2d}, d_{2d+1},
d_{2d+2}\}$; or
 \hfill\break
 (3): an induced 4-cycle $\{d_{2d}, d_{2d+l}, e_d,
e_{d+1}\}$.
 \hfill\break
 (Here, outer and inner wheel indices may be counted
modulo $2k$ and $k$, respectively.)
 \end{lemma}

 \begin{proof} The 3-stability essentially is a
reformulation of \cite[lemma~4.2.1]{RK}, and the
independence number is implicitly determined by
\cite[lemma~4.2.2~(b)]{RK}.

 For the proof of the rest, let $M$ be a
destabiliser of size 4 in \DW k. Since \DW k is
connected and contains at least eight trivalents, at
least one of them, say $v$, is a neighbour of $M$
(but not contained in $M$). Now, \DW k is a 3-stitch
of \w{k-1} with apex $v$, as depicted in figure IV,
p.~72, in \cite{RK}. Moreover, on the one hand,
there is some $x \in M \cap \lk v$, while on the
other hand $M' := M \cap V(\w{k-1})$ destabilises
\w{k-1}. (Else, there were an independent
$(k-1)$-subset $S$ of $V(\w{k-1}) \setminus M$,
whence $S \cup \{v\}$ were an independent $k$-subset
of $V(\DW k) \setminus M$; but $M$ destabilises \DW
k.)

 Thus and by lemma~\ref{L:A}, $M'$ consists of one
of the four bivalents in \w{k-1}, together with its
two neighbours in \w{k-1}, while $x$ is one of the
three neighbours of $v$. This leaves just twelve
potential $M$ to investigate, and it it easy to see
that most of them are not destabilisers of \DW k.
The lemma follows.
 \end{proof}

 \bigskip
 The next result should be rather well-known.
 \begin{lemma}
 \label{L:C}
 \5 is 4-stable.
 \end{lemma}

 \begin{proof}
 If $M \subset V(\5)$ and has at most four
vertices, then the induced graph on $V(\5)
\setminus M$ has at least $9 = R(3,4)$ vertices,
whence $M$ does not destabilise \5.
 \end{proof}

 \begin{lemma}
 \label{L:D}
 There is a 4-cycle through each vertex of degree at
least three in the \DW k ($k \ge 4$), \5, and the \w
k ($k \ge 2$) and other graphs in $\Gamma \cup
\Gamma'$. In particular, any graph in $\Gamma \cup
\Gamma'$ contains a 4-cycle, except $\0 \in \Gamma$
and $\c5 \in \Gamma'$.
 \end{lemma}

 \begin{proof}
 Inspection of the enumerated graphs.
 \end{proof}

 \bigskip \bigskip
 \section{The proof of the theorem.}
 \label{S:prf}
 As remarked in section~\ref{thm}, it is sufficient
to prove the following for each $k$, by means of
induction with respect to $k$:
 $$\hbox{If $\comp G = 1$, $\idn G = k$, $\nu \in
\{0,1\}$, and $G \in \Psi^{(\nu)}$, then } G \in
\Gamma^{(\nu)}. \leqno(6)$$
 Thus, for the whole proof, fix a positive integer
$k_0$, and assume that~(6) holds for each $k<k_0$,
that $G = (V,E) \in \Psi \cup \Psi'$, and that $k =
k_0 = \idn G$. Moreover, let $n = n(G)$, $e = e(G)$,
$t = t(G) \in \{0,1\}$, and $\mv = \mval G $ (where
$\mv\le2$ if $t=1$), and let $v$ be a vertex with
maximal second valency among the vertices with
minimal valency in $G$. In other words, we assume
that
 $$\v v=\mv \land \bigl( \v w=\mv \implies \vv2w \le
\vv2v \bigr) \,. \leqno (7)$$
 Finally, let the neighbours of $v$ be \row w{\v v},
where we may assume
 $$\v v \le \v{w_1} \le \ldots \le \v{w_{\v v}},
\leqno(8)$$
 and for $i = 1 \kdots \mv$, let $X_i = \lk(w_i)
\cap V(G_v)$, and let $X = \bigcup\limits_{i=1}^\mv
X_i = \B2v \setminus \B1v$.

 \smallskip
 If there were a redundant $\ep \in E$, then $(t=0
\implies t(V, E \setminus \{\ep\}) = -1)$, and $(t=1
\land \mv \le 2 \implies t(V, E \setminus \{\ep\}) =
0 \land \mval V,E\setminus\{\ep\} \le \mv \le 2)$,
in either case contradicting proposition~\ref{P:41}.
Thus, instead,
 $$G \hbox{ is edge critical.} \leqno (9)$$
 In particular, lemma~\ref{L:ec1um} applies, whence
 $$\idn{G_v} = k_0-1\,.$$

 Thus and by~(7), in particular, on the one hand
$\vv2v \ge \mv^2$, while on the other hand
 $$0 \le t(G_v) = \bigl( e-\vv2v \bigr) - 6
(n-\mv-1) + 13 (k_0-1) = t+6\mv-7 - \vv2v \,.$$
 Summing up, we have the useful restrictions
 $$\mv^2 \le \vv2v \le t+6\mv-7 \le 6\mv-6 \,.
\leqno (10)$$

 \medskip
 We start by considering a $G \in \Psi'$. Note that
then $\mv=2$, since a lower value would contradict
(10). For the same reason,
 $$4 \le \vv2v \le 6\,.$$
 We thus may make a case division with respect to
the value of \vv2v. However, first note that for
either value
 $$\comp{G_v}=1 \leqno (11)$$
 by~(9) and lemma~\ref{L:prpt}.

 \medskip
 $\underline{\vv2v = 2+2 = 4}$: By~(10) and~(7) both
$\v{w_i}=2$ and both $\vv2{w_i} = 4$, too. In other
words, each bivalent only has bivalent neighbours,
whence the connected graph $G$ must be 2-regular.
Thus, $G = \c l$ for some $l\ge4$. In fact, we must
have $l=5$, by the arguments in the proof of
\cite[lemma~2~(a)]{RK'} (or by directly calculating
the $t(\c l)$ for all $l$). Thus, indeed, $G \simeq
\c5 = \w2 \in \Gamma'$.

 \medskip
 $\underline{\vv2v = 2+3 = 5}$: By~(8), then
$\v{w_1} = 2$ and $\v{w_2} = 3$. Moreover,
$t(G_v) = 1$, too, whence (6) applies inductively
for $G_v$. In particular, thus $\mval G_v \ge 2$.
Hence, if $x$ is the single neighbour of $w_1$ in
$G_v\,$, then
 $$\v x \ge 1+ \vind {G_v}x \ge 1+2 = 3 \land
\vv2{w_1} = \v v+\v x \ge 2+3 = 5\,.$$
 On the other hand, $\vv2{w_1} \le \vv2v = 5$
by the assumption (7). We thus must have equalities.

 In particular, $x$ is a bivalent in $G_v$, which
thus belongs to $\Psi'$, and thus by~(6) to
$\Gamma'$. Thus and by~(11), up to isomorphisms,
either $G_v = \w{k_0-1}$, or
    $$G_v = \crt^{i+1}( \DW{k_0-2-i};M )$$
 for some $i \ge 0$ and a destabilising 4-subset $M
\subset V(\DW{k-2-i})$ of one of the three kinds
enumerated in lemma~\ref{L:B}. In either case, it is
sufficient to prove that $G$ is a 2-stitch with apex
$v$ and based at $x$ and its two neighbours in
$G_v$. In other words, we want to show that
 $$X = \B1{G_v,x},$$
 or, equivalently, that
 $$X_2 = \lk_{G_v}(x).$$

 Now, $X$ is a destabilising subset of $G_v$ of size
$\le3$, since the induced graph on $V(G_v) \setminus
X$ is $G_{w_1,w_2}$ and
 $$\idn {G_{w_1,w_2}} \le k_0-2 < \idn {G_v}$$
 by lemma~\ref{L:23}. Moreover, since $G$ is
triangle free, $X_2$ is an independent 2-set in
$G_v$. Thus, we have a trichotomy: Either $\ca X =
2$, or $\ca X = 3 \land \comp X \ge 2$, or indeed $X
= \B1{G_v,x}$; we have to prove that the first two
alternatives are impossible.

 For $G_v = \w{k-1}$, the destabilisers of size
$\le3$ are characterised in lemma~\ref{L:A}, and
they indeed must be of size~3 and connected. Thus,
assume instead that $G_v =
\crt^{i+1}(\DW{k-2-i};M)$, and, for a contradiction,
that either $X = X_2$ is an independent 2-set, or
$\ca X = 3$ but $\ca {E_X} \le 1$.

 Since $\idn{G_{w_1,w_2}} \le k_0-2$ (and by
proposition~\ref{P:41}),
 $$\begin{array}{l}
 0 \le t(G_{w_1,w_2}) \le t(G_v) - \Bigl( \sum_{y
\in X} \vind {G_v}y \ - \ \ca {E_X} \Bigr) + 6 \ca X
- 13\\
 {\mkern9.2mu} = 6 (\ca X-2) + \ca {E_X} - \sum_{y
\in X} \vind {G_v}y.
 \end{array}$$
 Since moreover $\mval G_v = 2$, and $G_v$ (and thus
$X$) contains at most two bivalents, $\ca X=2$ would
yield a glaring contradiction, and $\ca X=3$ but $X$
disconnected only could be possible if in addition
$\ca{E_X} = 1$ and $X$ consists of two bivalents and
one trivalent in $G_v$; and, moreover, then
$t(G_{w_1,w_2}) = 0$, whence $G_{w_1,w_2} \in
\Gamma$ by~(6)..

 However, then the single edge in $X$ must be the
edge between the two bivalents. Thus, if $y$ is the
trivalent in $X$, then all its three $G_v$
neighbours (say $z_1$, $z_2$, and $z_3$) belong to
$V(G_{w_1,w_2})$. Thus, either some $z_i$ were
trivalent in $G_v$, and therefore of valency less
than~3 in $G_{w_1,w_2}$, contradicting $\mval
G_{w_1,w_2} \ge 3$ (by proposition~\ref{P:41}); or
 $$\vvind {G_v}2y = \sum_i \vind {G_v}{z_i} \ge
3\cdot4 = 12 \implies t(G_{v,y}) \le 0,$$
 but since $\vind {G_{v,y}}x \le 2 < 3$, again we
would have a contradiction to the proposition.

 Thus, indeed we have eliminated all possibilities,
with the exception $X = \B1{G_v,x}$, whence indeed
$G = \crt(G_v;x) \in \Gamma'$.

 \medskip
 $\underline{\vv2v = 6}$: Then $t(G_v) = 0$, but
$G_v$ is destabilised by $X$. By~(11), (6), and
lemmata~\ref{L:C} and~\ref{L:B}, thus indeed $G_v =
\DW{k_0-1}$ and $G = \crt(\DW{k_0-1};X)$, with $X$
being of one of the three kinds given in
lemma~\ref{L:B}. Thus, indeed, then $G \in \Gamma'$.

 \bigskip
 Thus, we have proved~(6) for $G \in \Psi'$ (and $k
= k_0$). In particular, if there is a counterexample
to~(6) with a minimal value $k = k_0$ for its
independence number, then $G \in \Psi \setminus
\Gamma$, whence in particular
 $$k_0 \ge 7$$
 by proposition~\ref{P:41}. In the rest of the
proof, for a contradiction, we indeed assume that
$G$ is such a minimal counterexample.

 Since $t=0$, and by proposition~\ref{P:41}, and
summing up, we may strengthen (10) somewhat: In the
sequel we may assume that
 $$t=0 \land k_0 \ge 7 \land 3 \le \mv \le 4 \land
\mv^2 \le \vv2v \le 6\mv-7\,.$$
 This leaves five cases to consider, with $\mv=3$
and $\vv2v \in \{9,10,11\}$, and with $\mv=4$ and
$\vv2v \in \{16,17\}$. Again, we mainly consider
them separately. We start with the $\mv=3$ cases.

 \medskip
 $\underline{\vv2v=9}$: As in the $\vv2v=\mv^2=4$
case, this would force $G$ to be regular; this time,
3-regular, contradicting proposition~\ref{P:41}.

 \medskip
 Temporarily suspending the case division analysis,
note that in the remaining two $\mv=3$ cases,
 $$0 \le t(G_v) = t - \vv2v + 6(\v v+1) -13 =
11-\vv2v \le 1\,. \leqno (12)$$
 Furthermore, if in addition $G_v = G'+G''$ with the
$G^{(\nu)}$ non-empty, then without loss of
generality we could assume $t(G') = 0 \le t(G'') \le
11-\vv2v \le 1$.

 However, if then moreover $\vv2v = 11$, then both
$G^{(\nu)}$ were contained in $\Gamma$ by~(6), and
thus were strongly 3-stable, by lemmata~\ref{L:B}
and~\ref{L:C} for their components. Thus and by
lemma~\ref{L:prpt}~$(b)$, then $11 = \vv2v \ge 3+3+6
= 12$, a contradiction.

 Likewise, if then instead moreover $\vv2v = 10$,
then necessarily $\mval G'' = \mval G_v \le 2$, as
we shall see in a moment (and since $\mval G' \ge 3$
by proposition~\ref{P:41}), whence then $G'$ were
strongly 3-stable, and $G''$ would belong to
$\Gamma'$ and thus be strongly 2-stable, by also
employing lemma~\ref{L:A} for one $G''$ component;
this time yielding the contradiction $10 \ge 3+2+6 =
11$.

 Thus, instead, in the remaining $\mv=3$ cases, we
again have
 $$\comp{G_v}=1\,. \leqno (13)$$

 \smallskip
 $\underline{\vv2v=3+3+4=10}$: By~(8) and~(7), then
$\v{w_1} = \v{w_2} = 3$ $\v {w_3} = 4$, and moreover
$\vv2{w_1} \le 10 \ge \vv2{w_2}$, too. Thus, $X_1$
and $X_2$ contain trivalents, whence indeed $\mval
G_v \le 2$. Thus and by~(12), (13), and inductively
by~(6), then $G_v \in \Gamma'$ (and in fact $G_v$
contains at least two bivalents), and more precisely
$G_v = \w{k_0-1}$, or $G_v = \crt^{i+1} (\DW j;M)$
where $M \in \{K_{1,3}, \p4, \c4\}$ were a
destabiliser of \DW j, and $j = k_0-(i+2) \ge 5-i$.

 However, the latter possibility may be discarded:
If so, then no trivalent in $G_v$ could remain a
trivalent in $G$, since then some such remaining
trivalent $x$ would have $\vv2x > 10$; but there
were at least $2j - 4 + 2i \ge 6$ trivalents in
$G_v$; but $X$ would not contain more than seven
elements, including the two bivalents of $G_v$,
whence some $G_v$ trivalent indeed would be outside
$X$ and thus remain trivalent in $G$.

 Thus, instead, $G = \w{k_0-1}$. Hence,
\cite[lemma~4.2.2]{RK} yields that $G \simeq \DW k
\in \Gamma$, against the assumption that it were a
minimal counterexample.

 \medskip
 $\underline{\vv2v=11}$: Since $t(G_v) = 0$ and
by~(6) and~(13), and since $k_0 \ge 7 > 5$,
 $$G_v \in \Gamma \land G_v \not\simeq \5 \implies
G_v = \DW{k_0-1}\,.$$

 Next, note that the cardinality of $X$ is at most
8, and that equality holds if and only if the $X_i$
are disjoint, i.e., if and only if $\o{\c4}{G,v}=0$.
Moreover, any $w_i$ of degree 3 has $\vv2{w_i} \ge
3+4+4 = 11$, whence we have equality by~(7), and
thus then
 $$t(G_{w_i}) = 0 \implies G_{w_i} \in \Gamma
\implies \mval{G_{w_i}} \ge 3;$$
 whence there cannot be more than one trivalent
$w_i$. Thus and by~(8), and analogously,
 $$\v{w_1} = 3 \land \v{w_2} = \v{w_3} = 4 \land (x
\in X_1 \implies \v x=4).$$

 Next, the number of trivalents in $G_v =
\DW{k_0-1}$ is $2(k_0-1) \ge 12$. Thus, some of the
$G_v$ trivalents are adjacent to $X$ but not
contained therein. Each such trivalent $y$ must have
$\vv2y = 3+4+4$ in $G$; and they appear in pairs.
Reciprocially, for any such $y$, $G_y \simeq
\DW{k-1}$, too, and $v \in V_y$. Thus and by
lemma\ref{L:D},
 $$\o{\c4}{G,v} \ge \o{\c4}{G_y,v} \ge 1;$$
 whence $X$ has at most 7 vertices; whence there are
at least $k_0-4$ pairs of trivalents in $G_v$, which
remain trivalent in $G$. However, some of these
pairs would have to be `too close' in $G_v$,
yielding a trivalent $y$ with $t(G_y) = 0$ but
$\mval G_y < 3$, a contradiction.

 \bigskip
 Thus, instead, $\v v = \mv  = 4$, and $16 \le \vv2v
\le 17$. We treat the higher value first, since it
can be done briefly.

 \medskip
 $\underline{\vv2v = 4+4+4+5 = 17}$: By~(8),
$\v{w_1} = \v{w_2} = \v{w_3} = 4$, but $\v{w_4} =
5$. Moreover, since $t(G_{w_4}) \ge 0$ by
proposition~\ref{P:41}, $\vv2{w_4} \le 23$. However,
then $X_4$ must contain a tetravalent, say $u$; and
necessarily $\vv2u = 4+4+4+5 = 17$, too. This would
make $t(G_u) = 0$ and $v$ an element of $V(G_u)$,
whence by the inductive assumptions and
lemma~\ref{L:D}, there were a 4-cycle going through
$v$ and two of its neighbours of degree 4, say $w_1$
and $w_2$.

If $x$ were the last vertex of that 4-cycle, we
would be in a dilemma, as regards the degree of $x$.
Either, $x$ were tetravalent, and therefore of
degree at most 2 in $G_v$; or it were pentavalent,
whence $\vv2{w_1} = 17$ and $t(G_{w_1}) = 0$, but
$w_2$ were of degree at most 2 in $G_{w_1}$. In
either case, we would get a contradiction to
proposition~\ref{P:41}.

 \medskip
 $\underline{\vv2v = \mv^2 = 16}$: In analogy with
the other $\vv2v=\mv^2$ cases, $G$ is 4-regular.
Moreover, for each vertex $x$, $t(G_x) = 1$, whence
and inductively by~(6) either $\mval G_x > 2$ or
$G_x \in \Psi' \implies G_x \in \Gamma' \implies
\mval G_x = 2$. In other words, anyhow,
 $$\mval G_x \ge 2, \ \hbox{ and }\ \nv3{G_x} +
2\nv2{G_x} = \vv2x - \v x = 12\,. \leqno (14)$$

 It remains to show that indeed $G \simeq \5$. Now,
if in addition there is a 4-cycle in $G$, going
through $v$, say, then this is relatively easy to
show. As we just saw, then $G_v \in \Psi'$ and
$\mval G_v = 2$. Moreover, we in addition may assume
$v$ to be chosen in such a manner that the
independence number of the induced graph on the set
of $G_v$ bivalents is as large as possible.

 Now, since $G$ is connected and 4-regular, $G_v$
has no 4-regular component. However, for each not
4-regular connected graph $H$ in $\Gamma \cup
\Gamma'$, $\nv3H + 2\nv2H \ge 8$. Thus, and by~(14)
and linearity,
 $$ 12 = \nv3{G_v} + 2\nv2{G_v} \ge 8 \comp G\,;$$
 whence $G_v$ must be connected, and in fact,
without loss of generality, either $G_v = \w3$, or
$G_v = \crt^i(\DW{7-i};\{d_3,d_4,d_5,d_6\})$ for
some $i \in \{1,2,3\}$.

 However, in the latter case, $G_{d_4}$ would
contain an independent 2-set of bivalents,
consisting of $e_2$ and of one of the vertices added
to \DW{7-i} in the first stitch of the chain, but
$G_v$ would contain no such independent 2-set, in
contradiction to the choice of $v$.

 Thus, instead, in fact $G_v \simeq \w3$, from where
it is easy to deduce by inspection that indeed $G
\simeq \5$.

 \bigskip
 Thus, only the seemingly hardest case remains, that
$G$ would be both connected, 4-regular, and square
free. Actually, in \cite{RK}, the same kinds of
potential counterexamples, but to the statement (1),
also gave rise to considerable work; Radziszowsky
and Kreher use nine pages just to eliminate this
case\footnote{In \cite B, the elimination is
referred to a slightly more general result, whose
proof is even longer.}, consisting of their section
5.2, and of section~5.3 to the end of the proof of
their main theorem. Happily enough, most of their
proof also works in our situation, with the help of
a few observations.

 To be more precise, Radziszowski and Kreher
prove~(1) by induction, assuming that indeed $t(H)
\ge 0$ for any triangle free $H$ with $\idn H <
k_0$, and then consider a potential counterexample
 $$G \in \Lambda := \{G : \cqn G\le2 \land \idn
G=k_0 \land t(G)<0\}.$$
 They reasonably fast prove that then $G$ must be
connected and 4-regular, and must have $\girth G \ge
5$ (\cite[lemma 5.1.6 and proposition 5.1.8]{RK}).
They then eventually prove that the existence of
such a $G$ yields a contradiction. The short story
is that the proof of the latter {\sl mutatis
mutandis\/} may be applied for our $G \in \Gamma$.
Granting this, the induction step and thus our main
theorem is proved.

 \smallskip
 The somewhat longer story is that some care should
be taken with the ``things to be changed''. I
therefore provide a `translation' of their sequence
of lemmata to our situation, with the emphasis on
the changes, and omitting all parts of the proofs
which indeed are unchanged. In particular, I
introduce a shortcut, simplifying the treatment of
6-cycles. My intention is that my summary should be
intelligible in itself; however, a comprehensive
understanding of the full proof probably is hard
without accessing \cite{RK} directly.

 In fact, while the \cite{RK} arguments repeatedly
employ that their $G$ has $t(G) = -1$, this is
mainly used indirectly. They prove, that $n(G)$ must
be fairly large (\cite[proposition 5.1.7]{RK}), and
that thus certain subgraphs of the form $G_{\row
vr}$ both must be non-empty and have the $t(G_{\row
vr}) \ge t(G)+1$, with $\mval G_{\row vr} \ge 3$ in
case of equality. Moreover, in all the applications,
they {\sl a fortiori\/} are able to exclude
$\idn{G_{\row vr}}$ strictly less than the bound
given by lemma~\ref{L:23}.
 We start by proving three statements substituting
for this, and then show how to use them in order to
modify the lemma proofs in~\cite{RK}.

 Thus, again, let $G \in \Gamma$ be an assumed
minimal counterexample to theorem~\ref{T:1}, with
$\idn G = k_0$, $n(G) = n$, and $e(G) = e$, and
recall that then $G$ is 4-regular and connected, and
has $\o{\c4}G = \o{K_3}G = 0$, i.~e., has $\girth(G)
\ge 5$. In particular, $e-6n+13k_0 = t = t(G) = 0$
and $e = 2n$, whence $4n = 13k_0$, and $k_0$ is
divisible by~4. Since moreover $k_0 \ge 7$ by
proposition~\ref{P:41}, we actually must have $k_0
\ge8$, and thus get
 $$n \ge 26, \leqno (15)$$
 which should replace the calls to
\cite[lemma~5.1.7]{RK} in Radzisowski's and Kreher's
proofs.

 For any non-empty independent set $S = \{\row vr\}$
in $G$, $\idn {G_{\row vr}} \le k_0-r < k_0$ by
lemma~\ref{L:23}, whence the inductive assumptions
yield that $G_{\row vr} \in \Gamma$ if $t(G_{\row
vr}) \le 0$, and $G_{\row vr} \in \Gamma'$ if
$t(G_{\row vr})=1$ and $\mval G_{\row vr} \le 2$.
However, since $G$ contains no 4-cycle, neither does
$G_{\row vr}$ whence $(G_{\row vr} \in \Gamma
\implies G_{\row vr} = \0)$, and likewise $(G_{\row
vr} \in \Gamma' \implies G_{\row vr} \simeq \c5)$.
Together with (15), this yields that for such an
$S$
 $$t(G_{\row vr}) \ge 1 \hbox{ if } \ca{ \B1S} < 26,
\leqno (16)$$
 and
 $$t(G_{\row vr}) \ge 2 \hbox{ if } \ca{ \B1S} < 21
\hbox{ and } \mval G_{\row vr} \le 2 \,. \leqno
(17)$$

 On the other hand, since $t(G) = t = 0$ and by
lemmata~\ref{L:stsg} and~\ref{L:23},
 $$\begin{array}{rcl}
 t(G_{\row vr}) &=& e(G_{\row vr})-e - 6(n(_{\row
vr})-n) + 13 ( \idn{G_{\row vr}}-k_0 )\\
 &\le& 6 \ca{\B1S} - \v{\B1S} - 13r + \ca{E_{\B1S}};
 \end{array}$$
 and indeed the main application of (16) and~(17) is
to provide lower bounds for $\ca{E_{\B1S}}$ (which
Radziszowsky and Kreher call the number of edges in
the {\sl support\/} of $S$). Calls to (16) and
to~(17) should replace calls to
\cite[formula~(4)]{RK} and to
\cite[lemma~5.1.5]{RK}, respectively.

 \smallskip
 We now list the sequence of properties for $G$,
which leads to a contradiction. In most of them,
Radziszowski and Kreher consider a fixed vertex $v$,
let $H = G_v$, and let $J$ be the induced graph on
the set of trivalents in $H$.\footnote{Radziszowski
and Kreher in parts of their proofs change meanings
of $n$ and $e$; however, here we retain $n = n(G)$
and $e = e(G)$ consistently.} Thus, $V(J) = X$, and
by~$(14)$, $n(J) = 12$. Most of the properties
concern the graph structure of $J$. Let $C$ be an
arbitrary component of $J$.

 \cite[Lemma 5.2.2 (a) and (b)]{RK} state that for
any path $(v,t,u)$ of length~2 in $G$,
 $$\o{\c5,c_1,c_2,c_3}{G,v,t,u} \ge 1 \hbox{ and }
\o{\c5,c_1,c_2}{G,v,t} \ge 3\,. \leqno (18)$$
 The first claim is proved by noting that $\ca
{E_{\B1{\{v,u\}}}} \ge 9$ by applying
\cite[formula~(4)]{RK} to $G_{v,u}$; replace this by
applying (16).

 \cite[Lemma 5.2.3]{RK} states
 $$\girth J > 5,\
 \mval J > 1,\
 (s \in C \land \vind Js = 1 \implies C \simeq
\p2),\
 C \simeq \p2 \lor \mval C \ge 2\,.$$
 The only modification to be made of the proof
concerns the reason for the following fact (which we
shall reuse later):
 $$\hbox{If } s \in V(J) \hbox{ and } \vind Js = 1,
\hbox{ then } \mval H_s = \mval G_{v,s} \ge 3;
\leqno (19)$$
 apply (17) for $S = \{v,s\}$.

 \cite[Lemma 5.2.4]{RK} states
 $$\o{\c6}J = 0;$$
 if $(a,b,c,d,e,f)$ were a 6-cycle in $J$, then
apply (16) for $S = \{v,a,c,e\}$.

 \cite[Lemma 5.2.5]{RK} is somewhat technical; it
states that if $t \in V(J)$, $x,y \in \lk_J(t)$, and
$\vind Jx = \vind Jy = 2$, whence without loss of
generality $\lk_J(x) = \{t,x_1,x_2\}$ and $\lk_J(y)
= \{t,y_1,y_2\}$ with $x_2,y_2 \in V(H) \setminus
V(J)$, then
 $$x_1y_2, y_1x_2 \in E\,.$$
 The proof goes through without changes, as does the
proof of \cite[lemma~5.2.6]{RK}, stating
 $$\vind Jx = 3 \implies \vvind J2x = 2+2+3 = 7\,.$$

 These properties suffice to limit the possible $C$
to \p2, \c8, \c{10}, \c{12}, $S_1$, and $S_2$
(\cite[proposition~5.2.7]{RK}), where
 $$S_1 := (V(\c{12}), E(\c{12}) \cup \{c_6c_{12}\})
\hbox{ and } S_2 := (V(\c{12}), E(\c{12}) \cup
\{c_6c_{12}, c_3,c_9\}).$$
 In the next two lemmata, all of these except \p2
are discarded. The elimination of \c8 and \c{10}
goes through unmodified; for $C = J \in \{\c{12},
S_2\}$, Radziszowsky and Kreher prove that $n \le
20$, which here contradicts (15); and for $C = J =
S_1$, (16) should be applied for $S = \{v, c_1, c_3,
c_5, c_7, c_9, c_{11}\}$, yielding far too many
edges in \B1S.

 Thus, we know that (for any $v$)
 $$J \simeq 6\p2; \leqno (20)$$
 and in particular may deduce a sharper variant of
(18) (\cite[corollary~5.2.10 (a) and (b)]{RK}):
With $(v,t,u)$ as before,
 $$\o{\c5,c_1,c_2,c_3}{G,v,t,u} = 1 \hbox{ and }
\o{\c5,c_1,c_2}{G,v,t} = 3\,; \leqno (21)$$
 and immediately may deduce (\cite[corollary 5.2.11
(a)]{RK}
 $$\hbox{two 5-cycles can share at most one edge.}
\leqno (22)$$

 Radziszowsky and Kreher now proceed to investigate
6-cycles in $G$ in some detail. However, actually,
(19) and~(20) suffice to eliminate any such
6-cycle immediately: If instead \c6 were a subgraph
of $G$, then choosing $v := c_1$ and $s := c_3$, we
would have $s \in J$, hence $\vind Js = 1$
by~(20), and hence $\mval G_{v,s} \ge 3$
by~(19); but $\vind {G_{c_1,c_3}}{c_5} \le 2$, a
contradiction. Thus, indeed we have
 $$\o{\c6}G = 0 \,. \leqno (23)$$

 Thus, in order to achieve the final contradiction,
it is enough to prove that $G$ also must {\sl
contain \/} 6-cycles. Actually, in
\cite[lemma~5.3.2]{RK}, Radziszowsky and Kreher
proves that there would be at least six 6-cycles
through each edge $uv$ in their $G$. They start by
considering two 5-cycles $(u, v, x_1, \cdot, x_3)$
and $(u, v, x_2, \cdot, x_4)$ through $uv$ (existing
by~(21), and sharing no vertices outside $uv$
by~(22)), and then consider the independent set $S =
\{x_1, x_2, x_3, x_4\}$. Following their proof, but
applying (15) instead of
\cite[proposition~5.1.7]{RK} for $S$, we find that
also in our situation
 $$\ca {E_{\B1S}} \ge 21 > 19,$$
 indeed forcing the existence of 6-cycles, and thus
the sought contradiction to~(23).

 \medskip
 To sum up, we thus have proved, that if the claims
of theorem~\ref{thm} hold for all triangle free
graphs $G$ with $\idn G < k_0$, then they hold for
those with $\idn G = k_0$, too; whence indeed (3),
and thus the theorem, follows by induction.\noproof

 \bigskip \bigskip
 \section{Graphs $G$ with vertex numbers beyond
$\idn G$.}
 \label{S:disc}

 Finally, let us briefly discuss some possible
generalisations of \cite[theorem~5.1.1]{RK} and of
theorem~\ref{T:1}. In this survey section, some
proofs are omitted or just outlined; among these
those concerning the precise definition of graphs
who have crochet patterns with maximal valency at
most three. (However, all the local interpretations
of their patterns actually needed are presented in
section~\ref{S:step}, although in a somewhat
implicit manner.)

 For instance, it turns out, that a triangle free,
connected, and edge critical graph $G$ with $t(G)=1$
either has $\mval G = 2$ (and thus belongs to
$\Gamma'$), or has a crochet pattern $P$ with $\mval
P \ge 2$ and cycle space of dimension 2, or is one
of the two 4-regular $(3,6;16,32)$ realisers. (In
particular, hence, the last inequality in
corollary~\ref{cor} actually is strict.) On the
other hand, for $t(G)=2$, there are graphs
(including some $(3,9;26,52)$ realising ones), which
I do not know how to classify.

 \smallskip
 As we noted in section~\ref{S:bg}, (1) is just one
in a sequence of linear inequalities for triangle
free graphs, and the graphs with equalities in any
one of them are edge number critical. Conversely, in
the interval where equality may be attained in one
of the inequalities without violating the others,
all edge numbers critical graphs are of this type.
Thus, e.~g., the following statements are equivalent
for triangle free graphs $G$ with $2.5 \idn G \le
n(G) \le 3 \idn G$:

 \medskip\begingroup\sl
 \noindent$\begin{array}{ll}
 (i)&G \hbox{ is edge number critical with respect
to } \idn G \hbox{ and } n(G),\\
 (ii)&e(G)-5n(G)+10\idn G = 0,\\
 (iii)&G \hbox{ is a sum of chains and bicycles.}
 \end{array}$\endgroup

 \medskip
  (The equivalence between $(i)$ and $(ii)$ yields
the $5n-10k$ part of proposition~\ref{P:eval}.)

 \medskip
 Radziszowski and Kreher discussed whether there
could be an exhaustive set of intervals
 $a_i\idn G \le n(G) \le a_{i+1}\idn G$
 and linear inequalities
 $e(G)-b_in(G)+c_i\idn G \ge 0$,
 such that for each interval the edge critical
graphs are precise those with equality in the
corresponding inequality, at least for large enough
$n(G)$. Now, in this strong formulation, this
certainly is not the case. The `correct interval'
for the linear inequality~(1) would be
 $3\idn G \le n(G) \le 3.25 \idn G$,
 as seen from their own results and from
theorem~\ref{T:1}. However, if $3.25k-1 < n < 3.25k$
and $k \ge 5$, then $\e 3,k+1;n = 6n-13k+1$; and
such integers $n$ do exist for arbitrarily large $k$
not divisible by~4. Thus, the best we could hope for
is some kind of proportional result.

 With fixed $k$, $n$, and $e = \e3,k+1;n $, we may
compare the proportions $a = \frac nk$ and $b =
\frac ek$. By linearity, indeed then $\e3,mk+1;mn
\le me$ for any positive integer $m$. (In fact, if
$G$ is $(3,k+1;n,e)$ realising then clearly $mG$ is
$(3,mk+1;mn,me)$ realising.) By this and similar
considerations, in the limit we may consider $b$ as
a function of $a$, and this function is convex.

 Formally, for any real number $a\ge0$, let
 $$b(a) = \lim\limits_{k\rightarrow\infty} \frac{\e
{3,k+1;\mathopen\lfloor ak\mathclose\rfloor} }k$$
 (where as usual $\mathopen\lfloor
ak\mathclose\rfloor$ is the integer part of $ak$).
The enumerated linear inequalities imply that $b$ is
piecewise linear for $0 \le a \le 3.25$, and in fact
that there
 $$b(a) = \max(0, a-1, 3a-5, 5a-10, 6a-13) \,.$$
 What we could hope for, and what Radziszowski and
Kreher implicitly suggest, is, that $b$ continues to
be piecewise linear in its whole domain.

 To be more precise, they express their suggestions
in terms of the {\sl independence ratio\/} $\frac
kn$ as a function of the average degree
$\frac{2e}n$, instead of $\frac ek$ as a function of
$\frac nk$. Thus, in the limit, they define a
decreasing function $i^* = i^*(x)$. The relations
between the ratios are respected by the limits,
whence $i^*$ may be defined in terms of $b$, and
{\sl vice versa}. In fact, for $a>1$, if $x =
\frac{2b(a)}a$, then $i^*(x) = \frac1a$; and if $a =
\frac1{i^*(x)}$, then $b(a) = \frac x{2i^*(x)}$.
Thus, close to a point $(a_0,b_0)$ and the
corresponding point $(x_0,i^*_0)$, we have
 $$\begin{array}{rc}
 &i^* \hbox{ is linear in } x\\
 \iff&x \hbox{ is linear in } i^*\\
 \iff&\frac x{i^*} \hbox{ is linear in }
\frac1{i^*}\\
 \iff&b \hbox{ is linear in } a.
 \end{array}$$
 Thus, indeed the function $i^*$ is piecewise linear
(in its whole domain) if and only if $b$ is.

 However, whether or not these functions indeed are
piecewise linear seems to be a hard question.
 On the other hand, it is not very hard to prove
that $b$ is continous, convex, and non-decreasing,
and that it has both left and right derivatives in
its whole domain (except to the left at the origin),
and that moreover these derivatives are
non-decreasing, and that for any positive $a_0$ the
left derivative for $a_0$ is less than or equal to
the right derivative for $a_0$, but greater than or
equal to the right derivative for any $a < a_0$. In
fact, all the other properties are consequences of
the convexity and of $b$ being constant in an
interval containing the minimum of its domain.

 \begin{figure*}[t]
 \begin{picture}(200,360)(-80,0)

 \put(30,340){The function graph for $b(a)$, as far
as known:}


 \thinlines
 \put(60,10){\vector(1,0){180}}
 \put(60,10){\vector(0,1){300}}
 \put(245,7){$a$}
 \put(58,315){$b$}

 \put(99,0){$\scriptstyle1$}
 \put(138,0){$\scriptstyle2$}
 \put(155,0){$\scriptstyle2.\scriptscriptstyle5$}
 \put(177.7,0){$\scriptstyle3$}
 \put(185,0){$\scriptstyle3.\scriptscriptstyle25$}
 \put(218,0){$\scriptstyle4$}
 \put(50,47.5){$\scriptstyle1$}
 \put(50,87.5){$\scriptstyle2$}
 \put(45,107.5){$\scriptstyle2.\scriptscriptstyle5$}
 \put(50,127.5){$\scriptstyle3$}
 \put(50,167.5){$\scriptstyle4$}
 \put(50,207.5){$\scriptstyle5$}
 \put(50,247.5){$\scriptstyle6$}
 \put(45,267.5){$\scriptstyle6.\scriptscriptstyle5$}
 \put(50,287.5){$\scriptstyle7$}


 \thicklines
 \put(60,10){\line(1,0){40}}
 \put(100,10){\line(1,1){40}}
 \put(140,50){\line(1,3){20}}
 \put(160,110){\line(1,5){20}}
 \put(180,210){\line(1,6){10}}

 \put(100,10){\circle*2}
 \put(140,50){\circle*2}
 \put(160,110){\circle*2}
 \put(180,210){\circle*2}
 \put(190,270){\circle*2}

 \put(140,10){\circle*2}
 \put(160,10){\circle*2}
 \put(180,10){\circle*2}
 \put(190,10){\circle*2}
 \put(60,50){\circle*2}
 \put(60,110){\circle*2}
 \put(60,210){\circle*2}
 \put(60,270){\circle*2}

 \end{picture}
 \end{figure*}

 For the convexity, it is sufficient to note, that
for any $0 \le a_0 < a < a_1$ and any $\epsilon > 0$
 $$b(a) \le \frac{b(a_1)-b(a_0)}{a_1-a_0} a +
\frac{b(a_0)a_1-b(a_1)a_0}{a_1-a_0} + \epsilon,$$
 by considering graphs $G_0$ and $G_1$ with
$n(G_i)/\idn{G_i}$ and $e(G_i)/\idn{G_i}$
approximating $a_i$ and $b(a_i)$ sufficiently well,
and a suitable linear combination $G =
m_0G_0+m_1G_1$ with $n(G)/\idn G$ slightly larger
than $a$.

 That the graph of $b$ is convex also may be
reformulated thus:
 For every $c \in \mathopen[0,\infty\mathclose[\,$,
there is a unique number $d \in
\mathopen[0,\infty\mathclose[\,$, such that the line
$b-ca+d = 0$ touches but does not intersect the
graph. In other words,
 $$d = d(c) = \max \,(y : b(a) \ge ca-y
\;\forall\,a)
 = \min \,(y : \exists\,a \hbox{ such that } b(a) =
ca-y)\,.$$
 Hence, for each $c$ there is a ``best'' linear
graph invariant $\ell_c(G) = e(G)-cn(G)+d\alpha G$,
which is non-negative for all triangle free graphs.
(In fact, if there were a triangle free $G$ with
$\ell_c(G) < 0$, then $b(a) < ca-d$ for $a =
\frac{n(G)}{\idn G}$, against the assumptions.)
Correspondingly, with the same $c$ and $d$,
 $$\textstyle i^*(x) \ge \frac cd - \frac1{2d}x$$
 is an `optimal' linear lower bound for
Radziszowski's and Kreher's function $i^*$.

 For each such $c$, there either is a unique $a =
a_c$ with $b(a) = ca-d$, or there are several such
$a$. In either case, necessarily $c$ is at least the
left derivative and at most the right derivative of
$b$ at each such $a$. Thus, if there are several
such $a$ for a fixed $c$, they indeed form an
interval, and $b$ is differentiable with $b'(a) = c$
in the interior of this interval, and thus is linear
there.

 \smallskip
 Thus, in order to determine $b$ also for some
interval $3.25 \le a \le A$, we may equivalently
determine the left and right derivatives in this
interval (excepting the right derivative at $A$), or
determine the $\ell_c$ for all $c$ less than or
equal to the left $b$ derivative at $A$. The first
question would be the value $r$, say, of the right
derivative of $b$ at $a = 3.25$.
 Since this is at least equal to the left
derivative, and since on the other hand e.~g.\
$\e3,7;21 = 51 \implies b(3.5) \le 8.5$,
 $$6 \le r \le 8\,.$$
 In fact, for $3.25k \le n \le 3.5k$, there are
numerous realisers of $(3,k+1;n,8n-19.5k)$;
including all linear combinations of \5 and crochet
graphs with crochet patterns containing only
trivalent vertices. On the other hand, there is in
my knowledge no known triangle free graph $G$ for
which
 $$\ell(G) := e(G)-8n(G)+19.5\idn G$$
 is negative. Indeed, there is no such graph with
$\idn G < 10$, as can be seen from the exact values
and estimates of $e$-numbers in \cite{GR}. In fact,
also employing the $e$-number tables in \cite B, the
smallest possible counterexample would be a
$(3,11;35,84)$ realiser $G$; if it did exist, it
would have $\ell(G) = -1$. In view of the behaviour
for lower independence and vertex numbers, it would
be a bit of a surprise to have such a large simplest
example of a negative value for $\ell$.

 Thus, I think it is a reasonable guess that $r =
8$, or, equivalently, that $\ell_8 = \ell$, or
 $$b(a) \buildrel ?\over= 8a-19.5 \hbox{ for } 3.25
\le a \le 3.5\,.$$
 I have tried to prove this by means of the same
kind of strategy as the one used in this article,
but this seems hard. My best result so far is that
 $$r \ge 6.8\hbox{, i.e., } e(G) - 6.8e(G) + 15.6
\idn G = \ell_{6.8}(G) \ge 0$$
 for all triangle free $G$
 (\cite[proposition 12.5]B\footnote{Actually, in
order to work with integer valued invariants, the
proposition is formulated in terms of an invariant
$c(G) := 5\ell_{6.8}(G) = 5e(G)-34n(G)+78\idn G$.}
and its first corollary). However, both the
proposition and its proof are fairly complex; I
found no simpler way than making a simultaneous
induction over statements for graphs $G$ with over
forty distinct upper bounds on $\ell_{6.8}(G)$.
(Thus, just the formulation of the proposition
covers three typeset pages, and the proof of the
induction ten times more, preparatory results
uncounted.)

 Granted these results, we at least have e.~g.\ that
$8.2 \le b(3.5) \le 8.5$. Minor improvements of the
bounds for $r$, and for $b$ in the interval
$\mathopen[3.25, 3.5\mathclose]$ should be possible
with these methods; but for substantial
improvements, probably new ideas are needed, or an
improved interaction between theoretical analysis
and computer assisted investigation.

 \bibliographystyle{plain}

\begin{thebibliography}{9}
 \bibitem B
 J. Backelin, Contributions to a Ramsey
calculus, {\sl unpublished manuscript}.

 \bibitem{GR} J. Goedgebeur and S. P. Radziszowsky,
{\sl New computational upper bounds for Ramsey
numbers $R(3,k)$}, {\sl Electronic Journal of
Combinatorics}~20(1) (2013).

 \bibitem L A. Lesser, Theoretical and computational
aspects of Ramsey theory, {\sl Examensarbeten i
Matematik}, Matematiska Institutionen, Stockholms
Universitet~{\bf3} (2001).

 \bibitem{RK'} S. P. Radziszowski and D. L. Kreher,
On (3,k) Ramsey graphs: Theoretical and
computational results, {\sl J. Comb Math. and Comb.
Computing}~{\bf4} (1988), 37--52.

 \bibitem{RK} S. P. Radziszowski and D. L. Kreher,
Minimum triangle-free graphs, {\sl Ars
Comb}.~{\bf31} (1991), 65--92.
 \end{thebibliography}

 \end{document}